\documentclass[11pt,leqno]{article}
\usepackage{amsfonts,amssymb,amsmath,amsthm}
\usepackage{fullpage}

\newtheorem{theorem}{Theorem}[section]
\newtheorem{proposition}[theorem]{Proposition}
\newtheorem{lemma}[theorem]{Lemma}
\theoremstyle{definition}
\newtheorem{remark}[theorem]{Remark}
\numberwithin{equation}{section}

\newcommand{\N}{{\mathbb N}}
\newcommand{\R}{{\mathbb R}}
\newcommand{\grad}[1]{{\nabla\cdot{#1}\cdot\nabla}}
\newcommand{\supp}{\mathrm{Supp}}
\newcommand{\C}{{\mathcal C}}

\begin{document}

\title{Positive solutions to superlinear second--order\\
divergence type elliptic equations in cone--like domains}
\author{
{\Large Vladimir Kondratiev}\\
Department of Mathematics\\
and Mechanics\\
Moscow State University\\
Moscow 119 899, Russia\\
{\tt kondrat@vnmok.math.msu.su}\\
\and
{\Large Vitali Liskevich}\\
School of Mathematics\\
University of Bristol\\
Bristol BS8 1TW\\
United Kingdom\\
{\tt v.liskevich@bristol.ac.uk}\\
\smallskip\and
{\Large Vitaly Moroz}\\
School of Mathematics\\
University of Bristol\\
Bristol BS8 1TW\\
United Kingdom\\
{\tt v.moroz@bristol.ac.uk}}
\date{}
\maketitle

\begin{abstract}
We study the problem of the existence and nonexistence of positive solutions to
a superlinear second--order divergence type elliptic equation with measurable coefficients
$-\grad{a}u=u^p$ $(\ast)$, $p>1$, in an unbounded cone--like domain $G\subset\R^N$ $(N\ge 3)$.
We prove that the critical exponent
$p^*(a,G)=\inf\{p>1:\:\mbox{$(\ast)$ {\it has a positive supersolution at infinity in $G$}}\:\}$
for a nontrivial cone--like domain is always in $(1,\frac{N}{N-2})$ and in contrast
with exterior domains depends both on the geometry of the domain $G$ and the coefficients $a$ of the equation.
\end{abstract}

\section{Introduction}

We study the existence and nonexistence of positive solutions and supersolutions
to a superlinear second--order divergence type elliptic equation
\begin{equation}\label{*}
-\grad{a}u=u^p\quad\mbox{in }\:G.
\end{equation}
Here $p>1$, $G\subset\R^N$ ($N\ge 3$) is an unbounded domain
(i.e., connected open set) and
$-\grad{a}:=-\sum_{i,j=1}^{N}\frac{\partial}{\partial x_i}
\left(a_{ij}(x)\frac{\partial}{\partial x_j}\right)$
is a second order divergence type elliptic expression. We assume throughout
the paper that the matrix $a=(a_{ij}(x))_{i,j=1}^{N}$ is symmetric
measurable and uniformly elliptic, i.e. there exists an
ellipticity constant $\nu=\nu(a)>0$ such that
\begin{equation}\label{elliptic}
\nu^{-1}|\xi|^2\le\sum_{i,j=1}^{N}a_{ij}(x)\xi_i\xi_j\le\nu|\xi|^2,
\quad\mbox{for all $\xi\in\R^N$ and almost all $x\in G$.}
\end{equation}

The qualitative theory of semilinear equations of type (\ref{*})
in unbounded domains of different geometries has been extensively
studied because of various applications in mathematical physics
and reach mathematical structure.  One of the features of equation (\ref{*})
in unbounded domains is the nonexistence of positive solutions
for certain values of the exponent $p$.
Such nonexistence phenomena have been known at least since the celebrated
paper by Gidas and Spruck \cite{Gidas-Spruck}, where it was proved that
the equation
\begin{equation}\label{GS}
-\Delta u=u^p
\end{equation}
has no positive classical solutions in $\R^N$
($N\ge 3$) for $1\le p<\frac{N+2}{N-2}$. Though this results is sharp
(for $p\ge \frac{N+2}{N-2}$ there are classical positive solutions),
the critical exponent $p^\ast=\frac{N+2}{N-2}$ is highly unstable
with respect to any changes in the statement of the problem.
In particular, for any $p\in(\frac{N}{N-2},\frac{N+2}{N-2}]$
one can produce a smooth potential $W(x)$ squeezed between two positive constants
such that equation $-\Delta u=W(x)u^p$ has a positive solution in $\R^N$
(\cite{Toland}, see also \cite{Ding-Ni} for more delicate results).
If one looks for supersolutions to \eqref{GS} in $\R^N$ or studies \eqref{GS}
in exterior domains then the value and the properties of the critical exponent change.
The following result is well--known (see, e.g. \cite{Bandle-Levine,B-Veron}).
{\em If $N\ge 3$ and $1<p\le
\frac{N}{N-2}$ then there are no positive supersolutions to \eqref{GS}
outside a ball in $\R^N$. The value of the
critical exponent $p^\ast=\frac{N}{N-2}$ is sharp in the sense that
(\ref{*}) has (infinitely many) positive solutions outside a ball
for any $p>p^\ast$.}
This statement has been extended in different directions by many authors
(see, e.g. \cite{Bandle,BCN,B-Veron-Pohozaev,Birindelli,Deng-Levine,KLS,KLSU,
Levine,Pohozhaev,Serrin-Zou,Veron,Zhang-1,Zhang-2}).
In particular, in \cite{KLS} it was shown that the critical exponent
$p^\ast=\frac{N}{N-2}$ is stable with respect to the change of the Laplacian
by a second--order uniformly elliptic divergence type operator with measurable coefficients,
perturbed by a potential, for a sufficiently wide class of potentials
(see also \cite{KLSU} for equations of type (\ref{GS}) in exterior domains
in presence of first order terms).

In this paper we develop a new method of studying nonexistence of positive solutions
to \eqref{*} in cone--like domains (as model example of unbounded domains in $\R^N$
with nontrivial geometry). The method is based upon the maximum principle and asymptotic
properties at infinity of the corresponding solutions to the homogeneous linear equation.
This approach was first proposed in \cite{KLS}.
In the framework of our method we are able to establish the nonexistence results for \eqref{*}
with measurable coefficients in the cone--like domains without any smoothness of the boundary
in the setting of the most general definition of weak supersolutions.

We say that $u$ is a {\em solution} ({\em supersolution}) to equation (\ref{*}) if $u\in H^1_{loc}(G)$ and
$$\int_{G}\nabla u\cdot a\cdot\nabla\varphi\:dx=(\ge)
\int_{G}u^p\varphi\:dx\quad\mbox{for all }\: 0\le \varphi\in H^1_c(G),$$
where $H^1_c(G)$ stands for the set of compactly supported elements from $H^1_{loc}(G)$.
By the weak Harnack inequality for supersolutions (see, e.g. \cite[Theorem 8.18]{Gilbarg})
any nontrivial nonnegative supersolution to (\ref{*}) is positive in $G$.
We say that equation (\ref{*}) has a {\em solution (supersolution) at infinity}
if there exists a closed ball $\bar{B}_\rho$ centered at the origin with radius $\rho>0$
such that (\ref{*}) has a solution (supersolution) in $G\setminus \bar{B}_\rho$.

We define the {\em critical exponent} to equation (\ref{*}) by
$$p^\ast=p^*(a,G)=\inf\{p>1:\mbox{(\ref{*}) has a positive supersolution at infinity in $G$}\}.$$
In this paper we study the critical exponent $p^*(a,G)$
in a class of cone--like domains
$$\C_\Omega=\{(r,\omega)\in\R^N:\:\omega\in\Omega,\:r>0\},$$
where $(r,\omega)$ are the polar coordinates in $\R^N$
and $\Omega\subseteq S^{N-1}$ is a subdomain (a connected open subset)
of the unit sphere $S^{N-1}$ in $\R^N$.
The following proposition collects some properties of the critical exponent
and positive supersolutions to (\ref{*}) on cone--like domains.

\begin{proposition}\label{properties}
Let $\Omega^\prime\subset\Omega\subseteq S^{N-1}$ are subdomains of $S^{N-1}$. Then
\begin{enumerate}
\item[{$(i)$}]
$1\le p^\ast(a,\C_{\Omega^\prime})\le p^\ast(a,\C_\Omega)\le \frac{N}{N-2}$;
\item[{$(ii)$}] If
$p>p^\ast(a,\C_\Omega)$ then (\ref{*}) has a positive supersolution
at infinity in $\C_\Omega$;
\item[{$(iii)$}] If
$p>p^\ast(a,\C_\Omega)$ then (\ref{*}) has a positive solution
at infinity in $\C_\Omega$.
\end{enumerate}
\end{proposition}

\begin{remark}
Assertion (i) follows directly from the definition of the critical
exponent $p^\ast(a,G)$ and the fact that
$p^\ast(a,\R^N)=\frac{N}{N-2}$, see \cite{KLS}. Property (ii) simply
means that the critical exponent $p^\ast(a,G)$ divides the
semiaxes $(1,+\infty)$ into the nonexistence zone $(1,p^\ast)$ and
the existence zone $(p^\ast,+\infty)$. Existence (or nonexistence)
of a positive solution at the critical value $p^\ast$ itself is a
separate issue. Property (iii) says that the existence of a
positive supersolution at infinity implies the existence of a
positive solution at infinity. More precisely, we prove that if
(\ref{*}) has a supersolution $u>0$ in $\C_\Omega^\rho$ then for
any $r>\rho$ it has a solution $w>0$ in $\C_\Omega^r$ such that $w\le u$.
\end{remark}

The value of the critical exponent for the equation
$-\Delta u=u^p$ in $\C_\Omega$
with $\Omega\subseteq S^{N-1}$ satisfying mild regularity
assumptions was first established by Bandle and Levine \cite{Bandle-Levine} (see also \cite{Bandle}).
They reduce the problem to an ODE by averaging over $\Omega$.
The nonexistence of positive solutions without any smoothness
assumptions on $\Omega$ has been proved by Berestycki,
Capuzzo--Dolcetta and Nirenberg \cite{BCN} by means of a proper choice of a test function.

Let $\lambda_1=\lambda_1(\Omega)\ge 0$ be the principal eigenvalue of the
Dirichlet Laplace--Beltrami operator $-\Delta_\omega$ in $\Omega$.
Let $\alpha_-=\alpha_-(\Omega)<0$ be the negative root of the equation
$$\alpha(\alpha+N-2)=\lambda_1(\Omega).$$
The result in \cite{Bandle-Levine,BCN} reads as follows.
\begin{theorem}\label{Bandle}
Let $\Omega\subseteq S^{N-1}$ be a domain.
Then $p^\ast(id,\C_\Omega)=1-\frac{2}{\alpha_-}$,
and (\ref{*}) has no positive supersolutions at infinity in $\C_\Omega$
in the critical case $p=p^\ast(id,\C_\Omega)$.
\end{theorem}

Applicability of both ODE and test function techniques seems to be limited
to the case of radially symmetric matrices $a=a(|x|)$, whereas the method of the present paper
is suitable for studying equation (\ref{*}) with general uniformly elliptic
measurable matrix $a$.
It is extendable as far as the maximum principle is valid and appropriate asymptotic estimates
are available (see the proof of Theorem \ref{matrix} below).
Advantages of this approach are its transparency and flexibility.
As a first demonstration of the method we give a new proof of Theorem \ref{Bandle},
which has its own virtue being considerably less technical then in \cite{Bandle-Levine,BCN}.
As a consequence of Theorem \ref{Bandle} we derive the following result,
which says that in contrast to the case of
exterior domains the value of the critical exponent
on a fixed cone--like domain essentially depends on the coefficients
of the matrix $a$ of the equation.

\begin{theorem}\label{domain}
Let $\Omega\subset S^{N-1}$ be a domain such that $\lambda_1(\Omega)>0$.
Then for any $p\in(1,\frac{N}{N-2})$ there exists a uniformly elliptic matrix $a_p$
such that $p^\ast(a_p,\C_\Omega)=p$.
\end{theorem}

\begin{remark}
The matrix $a_p$ can be constructed in such a way that
(\ref{*}) either {\em has\,} or {\em has no} positive supersolutions at infinity in $\C_\Omega$
in the critical case $p=p^\ast(a_p,\C_\Omega)$,
see Remark \ref{R-log} for details.
\end{remark}

The main result of the paper asserts that equation (\ref{*}) with arbitrary
uniformly elliptic measurable matrix $a$ on a "nontrivial" cone--like domain
always admits a "nontrivial" critical exponent.

\begin{theorem}\label{matrix}
Let $\Omega\subseteq S^{N-1}$ be a domain and $a$ be a uniformly elliptic matrix.
Then $p^\ast(a,\C_\Omega)>1$.
If the interior of $S^{N-1}\setminus\Omega$ is nonempty
then $p^\ast(a,\C_\Omega)<\frac{N}{N-2}$.
\end{theorem}

The rest of the paper is organized as follows.
In Section 2 we discuss the maximum and comparison principles in a form appropriate
for our purposes and study some properties of linear equations in cone--like domains.
Proposition \ref{properties} is proved in Section 3.
Section 4 contains the proof of Theorems \ref{Bandle} and \ref{domain}.
The proof of Theorem \ref{matrix} as well as some further remarks are given in Section 5.

\section{Background, framework and auxiliary facts}

Let $G\subseteq\R^N$ be a domain in $\R^N$.
Throughout the paper we assume that $N\ge 3$.
We write $G^\prime\Subset G$ if $G^\prime$ is a subdomain of $G$ such that $cl\,G^\prime\subset G$.
By $\|\cdot\|_p$ we denote the standard norm in the Lebesgue space $L^p$.
By $c,c_1,\dots$ we denote various positive constants whose exact value is irrelevant.

Let $S^{N-1}=\{x\in\R^N:|x|=1\}$ and $\Omega\subseteq S^{N-1}$ be a subdomain of $S^{N-1}$.
Here and thereafter, for $0\le\rho<R\le+\infty$, we denote
$$\C_{\Omega}^{(\rho,R)}:=\{(r,\omega)\in\R^N:\:\omega\in\Omega,\:r\in(\rho,R)\},
\qquad\C_{\Omega}^\rho:=\C_{\Omega}^{(\rho,+\infty)}.$$
Accordingly, $\C_\Omega=\C_\Omega^0$ and $\C_{S^{N-1}}=\R^N\setminus\{0\}$.

\paragraph{Maximum and comparison principles}

Consider the linear equation
\begin{equation}\label{linear-f}
-\grad{a}u-Vu=f\quad\mbox{in }\:G,
\end{equation}
where $f\in H^1_{loc}(G)$ and $0\le V\in L^1_{loc}(G)$ is a {\em form--bounded potential}, that is
\begin{equation}\label{form-bound}
\int_G Vu^2\:dx\le
(1-\epsilon)\int_{G}\nabla u\cdot a\cdot\nabla u\:dx
\quad\mbox{for all }\: 0\le u\in H^1_c(G)
\end{equation}
with some $\epsilon\in(0,1)$.
A solution (supersolution) to (\ref{linear-f}) is a function $u\in H^1_{loc}(G)$ such that
$$\int_{G}\nabla u\cdot a\cdot\nabla\varphi\:dx-\int_G Vu\varphi\:dx=(\ge)
\langle f,\varphi\rangle\quad\mbox{for all }\: 0\le \varphi\in H^1_c(G),$$
where $\langle\cdot,\cdot\rangle$ denotes the duality between
$H^{-1}_{loc}(G)$ and $H^1_{c}(G)$.
If $u\ge 0$ is a supersolution to
\begin{equation}\label{linear-0}
-\grad{a}u-Vu=0\quad\mbox{in }\:G,
\end{equation}
then $u$ is a supersolution to $-\grad{a}u=0$ in $G$.
Therefore $u$ satisfies on any subdomain $G^\prime\Subset G$ the {\it weak Harnack inequality}
$$\inf_{G^\prime}u\ge\frac{C_W}{\mathrm{mes}(G^\prime)}\int_{G^\prime}u\:dx,$$
where $C_W=C_W(G,G^\prime)>0$.
In particular, every nontrivial supersolution $u\ge 0$ to (\ref{linear-0})
is strictly positive, that is $u>0$ in $G$.

We define the space $D^1_0(G)$ as a completion of $C^\infty_c(G)$
with respect to the norm $\|u\|_{D^1_0(G)}:=\|\nabla u\|_2$.
The space $D^1_0(G)$ is a Hilbert and Dirichlet space, with the dual $D^{-1}(G)$,
see, e.g. \cite{Fuk}. This implies, amongst other things, that $D^1_0(G)$ is invariant
under the standard truncations, e.g. $v\in D^1_0(G)$ implies that
$v^-=v\vee 0\in D^1_0(G)$, $v^+=v\wedge 0\in D^1_0(G)$.
By the Sobolev inequality
$D^1_0(G)\subset L^{\frac{2N}{N-2}}(G)$.
The Hardy inequality
\begin{equation}\label{Hardy}
\int_{\R^N}|\nabla u|^2\:dx\ge \frac{(N-2)^2}{4}\int_{\R^N}\frac{|u|^2}{|x|^2}\:dx
\quad\mbox{for all $u\in H^1_c(\R^N)$},
\end{equation}
implies that $D^1_0(G)\subset L^2(G,|x|^{-2}dx)$.
Since the matrix $a$ is uniformly elliptic and
the potential $V$ is form bounded, the quadratic form
$$Q(u):=\int_{G}\nabla u\cdot a\cdot\nabla u\:dx-\int_G Vu^2\:dx$$
defines an equivalent norm $\sqrt{Q(u)}$ on $D^1_0(G)$.
The following lemma is a standard consequence of the Lax--Milgram Theorem.

\begin{lemma}\label{existence}
Let $f\in D^{-1}(G)$. Then the problem
$$-\grad{a}v-Vv=f,\qquad v\in D^1_0(G),$$
has a unique solution.
\end{lemma}

The following two lemmas provide the maximum and comparison principles
for equation (\ref{linear-f}), in a form suitable for our framework.
We give the full proofs for completeness, though the arguments are mostly standard.

\begin{lemma}\label{maximum}
{\sf (Weak Maximum Principle)}
Let $v\in H^1_{loc}(G)$ be a supersolution to equation (\ref{linear-0})
such that $v^-\in D^1_0(G)$. Then $v\ge 0$ in $G$.
\end{lemma}

\proof
Let $(\varphi_n)\subset C^\infty_c(G)$ be a sequence such
that $\|\nabla(v^--\varphi_n)\|_2^2\to 0$. Set
$v_n:=0\vee\varphi_n\wedge v^-$. Since $0\le v_n\le v^-\in
D^1_0(G)$ and
\begin{eqnarray*}
\int_{G}|\nabla(v^--v_n)|^2\:dx&=&
\int_{\{0\le\varphi_n\le v^-\}}|\nabla(v^--\varphi_n)|^2\:dx+
\int_{\{\varphi_n\le 0\}}|\nabla v^-|^2\:dx\\
&\le&
\int_{G}|\nabla(v^--\varphi_n)|^2\:dx+ \int_{\{\varphi_n\le 0\}}|\nabla v^-|^2\:dx\to 0,
\end{eqnarray*}
by the Lebesgue dominated convergence, we conclude that $\|\nabla(v^--v_n)\|_2^2\to
0$ (cf. \cite[Lemma 2.3.4]{Fuk}). Taking $(v_n)$ as a sequence of
test functions we obtain
\begin{eqnarray*}
0&\le&
\int_{G}\nabla v\cdot a\cdot\nabla v_n\:dx-\int_G Vvv_n\:dx\\
&=&-\int_{G}\nabla v^-\cdot a\cdot\nabla v_n\:dx+\int_G Vv^-v_n\:dx\:
\to\:\int_{G}\nabla v^-\cdot a\cdot\nabla v^-\:dx+\int_G V|v^-|^2\:dx\le 0.
\end{eqnarray*}
Thus we conclude that $v^-=0$.
\qed

\begin{lemma}\label{comparison}
{\sf (Weak Comparison Principle)}
Let $0\le u\in H^1_{loc}(G)$, $v\in D^1_0(G)$ and
$$-\grad{a}(u-v)-V(u-v)\ge 0\quad\mbox{in }\:G.$$
Then $u\ge v$ in $G$.
\end{lemma}

\proof
Let $(G_n)_{n\in\N}$ be an exhaustion of $G$,
i.e.\ an increasing sequence of bounded smooth domains
such that $G_n\Subset G_{n+1}\Subset G$ and $\cup_{n\in\N}G_n=G$.
Let $v\in D^1_0(G)$. Let $f\in D^{-1}(G)$ be defined by duality as
$$f:=-\grad{a}v-Vv.$$
Let $v_n\in D^1_0(G_n)$ be the unique weak solution to the linear problem
$$-\grad{a}v_n-Vv_n = f,\qquad v_n\in D^1_0(G_n).$$
Then
$$-\grad{a}(u-v_n)-V(u-v_n)\ge 0\quad\mbox{in }\:G_n,$$
with
$$u-v_n\in H^1(G_n),\qquad 0\le(u-v_n)^-\le v_n^+\in D^1_0(G_n).$$
Therefore $(u-v_n)^-\in D^1_0(G_n)$.
By Lemma \ref{maximum} we conclude that $(u-v_n)^-=0$, that is $v_n\le u$.
Let $\bar{v}_n\in D^1_0(G)$ be defined as $\bar{v}_n=v_n$ on $G_n$,
$\bar{v}_n=0$ on $G\setminus G_n$.
To complete the proof of the lemma it suffices to show that
$\bar{v}_n\to v$ in $D^1_0(G)$.
Indeed,
$$Q(\bar{v}_n)=\int_{G}\nabla\bar{v}_n\cdot a\cdot\nabla\bar{v}_n-\int_{G}V|\bar{v}_n|^2 =
\langle f, v_n\rangle\le c\|f\|_{D^{-1}(G)}\|\bar{v}_n\|_{D^1_0(G)},$$
where $\langle\cdot,\cdot\rangle$ stands for the duality between $D^1_0(G)$ and $D^{-1}(G)$.
Hence the sequence $(\bar{v}_n)$ is bounded in $D^1_0(G)$.
Thus we can extract a subsequence, which we still denote by $(\bar{v}_n)$,
that converges weakly to $v_\ast\in D^1_0(G)$.
Now let $\varphi\in H^1_c(G)$.
Then for all $n\in\N$ large enough, we have that $Supp(\varphi)\subset G_n$ and
$$\int_{G}\nabla\bar{v}_n\cdot a\cdot\nabla\varphi-\int_{G}V\bar{v}_n\varphi =
\int_{G_n}\nabla v_n\cdot a\cdot\nabla\varphi-\int_{G_n}Vv_n\varphi=\langle f,\varphi\rangle.$$
By the weak continuity we conclude that
$$\int_{G}\nabla v_\ast\cdot a\cdot\nabla\varphi-\int_{G}Vv_\ast\varphi =
\langle f,\varphi\rangle.$$
Therefore $v_\ast\in D^1_0(G)$ satisfies
$$-\grad{a}v - Vv = f,\quad v\in D^1_0(\C_\Omega).$$
Hence $v_\ast=v$.
Furthermore,
$$Q(\bar{v}_n-v)=\langle f,\bar{v}_n\rangle-2\langle f,v\rangle+\langle f,v\rangle.$$
Since $\langle f,\bar{v}_n\rangle\to\langle f,v\rangle$
it follows that $\bar{v}_n\to v$ in $D^1_0(G_n)$.
\qed

\paragraph{Minimal positive solution in cone like domains}
Here we construct a minimal positive solution suitable for the framework of cone--like domains
(cf. Agmon \cite{Agmon}).
Let $\Omega\subseteq S^{N-1}$ be a domain. Consider the equation
\begin{equation}\label{min}
-\grad{a}u-Vu=0\quad\mbox{in }\:\C_\Omega,
\end{equation}
where $0\le V\in L^1_{loc}(\C_\Omega^\rho)$ is a form--bounded potential.
Let $0\le\psi\in C^\infty_c(\Omega)$
and $\theta_\rho\in C^\infty[\rho,+\infty)$
be such that $\theta(\rho)=1$, $0\le\theta\le 1$ and $\theta=0$ for $r\ge\rho+\epsilon$
with some $\epsilon>0$.
Thus $f_\psi:=\grad{a}(\psi\theta_\rho)\in D^{-1}(\C_\Omega^\rho)$.
Let $w_\psi$ be the unique solution to the problem
\begin{equation}\label{psi-tilde}
-\grad{a}w-Vw=f_\psi,\qquad w\in D^1_0(\C_\Omega^\rho),
\end{equation}
which is given by Lemma \ref{existence}.
Set $v_\psi:=w+\psi\theta_\rho$. Then $v_\psi$ is the solution to the problem
\begin{equation}\label{e-psi}
-\grad{a}v-Vv = 0,\qquad v-\psi\theta_\rho\in D^1_0(\C_\Omega^\rho).
\end{equation}
By the weak Harnack inequality $v_\psi>0$ in $\C_\Omega^\rho$.
We call such $v_\psi$ a
{\em minimal positive solution to equation (\ref{min}) in $\C_\Omega^\rho$}.
Notice that $v_\psi$ actually does not depend on the particular choice of the function $\theta_\rho$.
The minimality property of the above $v_\psi$ is shown in the next lemma.

\begin{lemma}\label{minimal}
{\sf (Minimal Solution Lemma)}
Let $v_\psi$ be a minimal positive solution to equation (\ref{min}) in $\C_\Omega^\rho$.
Then for any positive supersolution $u>0$ to equation (\ref{min}) in $\C_\Omega^r$
with $r\in(0,\rho)$ there exist $c>0$ such that
$$u\ge cv_\psi\quad\mbox{in }\:\C_\Omega^\rho.$$
\end{lemma}

\proof
Let $v_\psi:=w_\psi+\psi\theta_\rho$ be a minimal positive solution to (\ref{min}) in $\C_\Omega^\rho$.
Choose $\Omega^\prime\Subset\Omega$ such that $\supp(\psi)\Subset\Omega^\prime$.
Let $\epsilon>0$ be such that $\theta_\rho=0$ for all $r\ge\rho+\epsilon$.

Let $u>0$ be a positive supersolution to (\ref{min}) in $\C_\Omega^r$ with $r\in(0,\rho)$.
By the weak Harnack inequality there exists $m=m(\Omega^\prime,\epsilon)>0$ such that
$$u>m\quad\mbox{in }\:\C_{\Omega^\prime}^{(\rho,\rho+\epsilon)}.$$
Choose $c>0$ such that $c\psi<m$.
Then $u-c\psi\theta_\rho\ge 0$ in $\C_\Omega^\rho$, $cw_\psi\in D^1_0(\C_\Omega^\rho)$ and
$$(-\grad{a}-V)((u-c\psi\theta_\rho)-cw_\psi)=(-\grad{a}-V)u\ge 0\quad\mbox{in }\:\C_\Omega^\rho.$$
By Lemma \ref{comparison} we conclude that $u-c\psi\theta_\rho\ge cw_\psi$,
that is $u\ge cv_\psi$ in $\C_\Omega^\rho$.
\qed

\begin{remark}
Let $\Gamma_a(x,y)$ be the fundamental solution to the equation $-\grad{a}u=0$ in $\R^N$.
Then for any domain $\Omega\subseteq S^{N-1}$ the function $\Gamma_a(x,0)$ is a positive solution to
\begin{equation}\label{e-green}
-\grad{a}u=0\quad\mbox{in }\:\C_\Omega.
\end{equation}
By Lemma \ref{comparison} and the classical estimate \cite{LSW} we conclude that
any minimal positive solution $v_\psi$ to (\ref{e-green}) in $\C_\Omega^\rho$ obeys the upper bound
\begin{equation}\label{upper-green}
v_\psi\le c_1\Gamma_a(x,0)\le c_2|x|^{2-N}\quad\mbox{in }\:\C_\Omega^{\rho}.
\end{equation}
\end{remark}

\paragraph{Nonexistence Lemma}
The next lemma (compare \cite{KLS},\cite[p.156]{Pinsky})
is the key tool in our proofs of nonexistence of positive solutions
to nonlinear equation (\ref{*}).

\begin{lemma}\label{lin-nonexist}
{\sf (Nonexistence Lemma)}
Let $0\le V\in L^1_{loc}(\C_\Omega^\rho)$ satisfy
\begin{equation}\label{big}
|x|^2 V(x)\to\infty\quad
\mbox{as $x\in\C_{\Omega^\prime}^\rho$ and $|x|\to\infty$}
\end{equation}
for a subdomain $\Omega^\prime\subseteq\Omega$.
Then the equation
\begin{equation}\label{V}
-\grad{a}u-Vu=0\quad\mbox{in }\:\C_\Omega^\rho
\end{equation}
has no nontrivial nonnegative supersolutions.
\end{lemma}

The proof of this lemma is based upon the following simple result.

\begin{lemma}\label{BigPotential}
Let $G\subset\R^N$ be a bounded domain and $\lambda_1=\lambda_1(G)>0$ be
the principal Dirichlet eigenvalue of $-\grad{a}$ in $G$.
If $\mu>\lambda_1$ then the equation
\begin{equation}\label{lambda}
-\grad{a}u=\mu u\quad\mbox{in }\:G
\end{equation}
has no positive supersolutions.
\end{lemma}

\proof
Assume that $u\ge 0$ is a supersolution to (\ref{lambda}).
By the monotonicity property of the principal Dirichlet eigenvalue of $-\grad{a}$,
for any given $\mu>\lambda_1(G)$ one can choose a subdomain
$G^\prime$ such that $\lambda_1(G^\prime)<\mu$.
Let $e^\prime$ be the eigenfunction of $-\grad{a}$ corresponding to $\lambda_1(G^\prime)$,
$0<e_1^\prime\in H^1_0(G^\prime)\subset H^1_c(G)$.
Taking $e^\prime$ as a test function we have
$$0\le\mu\int_{G^\prime}u e_1^\prime\:dx\le
\int_{G^\prime}\nabla u\cdot a\cdot\nabla e_1^\prime\:dx=
\lambda_1(G^\prime)\int_{G^\prime}u e_1^\prime\:dx.$$
We conclude that $u=0$ in $G^\prime$ and,
by the weak Harnack inequality in $G$.
\qed

\paragraph{Proof of Lemma \ref{lin-nonexist}}
Let $\lambda_1(\C_\Omega^{(\rho,2\rho)})>0$ be the principal Dirichlet eigenvalue
of $-\grad{a}$ on $\C_\Omega^{(\rho,2\rho)}$.
Rescaling the equation
$-\grad{a}v=\lambda v$ from $\C_\Omega^{(\rho,2\rho)}$ to $\C_\Omega^{(1,2)}$
one sees that
$$\frac{c^{-1}}{\rho^2}\lambda_1(\C_\Omega^{(1,2)})\le
\lambda_1(\C_\Omega^{(\rho,2\rho)})\le
\frac{c}{\rho^2}\lambda_1(\C_\Omega^{(1,2)}),$$
where $c=c(a)>0$ depends on the ellipticity constant of the matrix $a$
and does not depend on $\rho>0$.

Let $u\ge 0$ be a supersolution to (\ref{V}).
Then (\ref{big}) implies that for some $R\gg 1$ one can find $\mu>0$ such that
$V(x)\ge\mu\ge c\lambda_1(\C_{\Omega^\prime}^{(1,2)})R^{-2}$ in $\C_{\Omega^\prime}^{(R,2R)}$.
Hence $u$ is a supersolution to
$$-\grad{a}u=\mu u\quad\mbox{in }\:\C_{\Omega^\prime}^{(R,2R)}$$
with $\mu>\lambda_1(\C_{\Omega^\prime}^{(R,2R)})$.
By Lemma \ref{BigPotential} we conclude that $u=0$ in $\C_{\Omega^\prime}^{(R,2R)}$.
Therefore by the weak Harnack inequality $u=0$ in $\C_\Omega^\rho$.
\qed

\section{Proof of Proposition \ref{properties}}

Property (i) is obvious. We need to prove (ii) and (iii).
\smallskip

$(ii)$
Let $p_0\ge p^\ast(a,\C_\Omega^\rho)$ be such that equation (\ref{*}) with exponent $p_0$
has a positive supersolution $u>0$ in $\C_\Omega^\rho$.
Let $p>p_0$ and $\alpha=\frac{p-1}{p_0-1}>1$.
Set $v:=u^{1/\alpha}$.
By the weak Harnack inequality $u>0$ in $\C_\Omega^\rho$.
Hence $u^{-s}\in L^\infty_{loc}(\C_\Omega^\rho)$ for any $s>0$.
Therefore
$\nabla v =\alpha^{-1}u^{1/\alpha-1}\nabla u\in L^2_{loc}(\C_\Omega^\rho)$,
that is $v\in H^1_{loc}(\C_\Omega^\rho)$.

Let $0\le \varphi\in C^\infty_c(\C_\Omega^\rho)$. Then
\begin{eqnarray*}
\int_{\C_\Omega^\rho}\nabla v^\alpha\cdot a\cdot\nabla\varphi\:dx&=&
\alpha\int_{\C_\Omega^\rho}v^{\alpha-1}\nabla v\cdot a\cdot\nabla\varphi\:dx\\
&=&\alpha\int_{\C_\Omega^\rho}\nabla v\cdot a\cdot\nabla(v^{\alpha-1}\varphi)\:dx-
\alpha(\alpha-1)\int_{\C_\Omega^\rho}\nabla v\cdot a\cdot\nabla v \:(v^{\alpha-2}\varphi)\:dx\\
&\le&\alpha\int_{\C_\Omega^\rho}\nabla v\cdot a\cdot\nabla(v^{\alpha-1}\varphi)\:dx.
\end{eqnarray*}
Notice, that $v^{\alpha-1}=u^{1-1/\alpha}\in H^1_{loc}(\C_\Omega^\rho)$ by the
same argument as above.
Therefore $v^{\alpha-1}\varphi\in H^1_c(\C_\Omega^\rho)$.
We shall prove that the set
$${\mathcal K}_v=\{v^{\alpha-1}\varphi,\:0\le\varphi\in H^1_c(\C_\Omega^\rho)\}$$
is dense in the cone of nonnegative functions in $H^1_c(\C_\Omega^\rho)$.
Indeed, let $0\le \psi\in H^1_c(\C_\Omega^\rho)$.
Let $\psi_n\in C^\infty_0(\C_\Omega^\rho)$ be an approximating sequence such that
$\|\nabla(\psi_n-\psi)\|_2\to 0$.
Set $\varphi_n=v^{1-\alpha}\psi_n^+$.
It is clear that $0\le \varphi_n\in {\mathcal K}_v\subset H^1_c(\C_\Omega^\rho)$
and $\|\nabla(v^{\alpha-1}\varphi_n-\psi)\|_2\to 0$.

Since $v^\alpha=u$ and $u>0$ is a supersolution of (\ref{*}),
we obtain that
$$\alpha\int_{\C_\Omega^\rho}\nabla v\cdot a\cdot\nabla(v^{\alpha-1}\varphi)\:dx\ge
\int_{\C_\Omega^\rho} v^{\alpha p_0}\varphi\:dx
=\int_{\C_\Omega^\rho} v^p(v^{\alpha-1}\varphi)\:dx$$
for any $0\le \varphi\in H^1_c(\C_\Omega^\rho)$.
Thus $\alpha^{1/(1-p)}v$ is a supersolution to equation (\ref{*}) in $\C_\Omega^\rho$
with exponent $p>p_0$.
\smallskip

$(iii)$
The existence of a (bounded) positive solution to equation (\ref{*}) with $p>\frac{N}{N-2}$ in $B_r^c$ for any $r>0$
has been proved in \cite{KLS}. We shall consider the case $p\le\frac{N}{N-2}$.

Let $u>0$ be a supersolution to (\ref{*}) with exponent $p\le N/(N-2)$ in $\C_\Omega^\rho$.
Fix $\psi\in C_c^\infty(\Omega)$ and $r>\rho$.
Let $v_\psi>0$ be a minimal solution in $\C_\Omega^r$ to $-\grad{a}v=0$ in $\C_\Omega$.
Then $u\ge cv_\psi$ in $\C_\Omega^r$ by Lemma \ref{comparison}.
Without loss of generality we assume that $c=1$.
Thus $v_\psi>0$ is a subsolution to (\ref{*}) in $\C_\Omega^r$ and $v_\psi\le u$ in $\C_\Omega^r$.
We are going to show that (\ref{*}) has a positive solution $w$ in $\C_\Omega^r$ such that
$v_\psi\le w\le u$ in $\C_\Omega^r$.

Let $(G_n)_{n\in\N}$ be an exhaustion of $\C_\Omega^r$.
Consider the boundary value problem
\begin{equation}\label{bvp}
\left\{\begin{array}{rcll}
-\grad{a}w&=&w^p&\mbox{in }\:G_n,\\
w&=&v_\psi&\mbox{on }\:\partial G_n.\\
\end{array}\right.
\end{equation}
Since $G_n\Subset\C_\Omega^r$ is a smooth bounded domain
and $v_\psi\in C^{0,\gamma}_{loc}(\C_\Omega^r)$, the problem (\ref{bvp}) is well--posed.
Clearly, $v_\psi\le u$ is still a pair of sub and supersolutions for (\ref{bvp}).
Notice that we do not assume that $u\in H^1(G_n)$ is bounded.
However, since $p\le \frac{N}{N-2}<\frac{N+2}{N-2}$,
one can use an $H^1$--version of sub and supersolution method, see e.g. \cite[Theorem 2.2]{Dancer}.
Thus there exists a weak solution $w_n\in H^1(G_n)$ of (\ref{bvp}) such that $v_\psi\le w_n\le u$ in $G_n$.

Consider a sequence $(w_n)_{n>1}$ in $G_1$.
Choose a function $\theta\in C^\infty_c(G_2)$ such that $0\le\theta\le 1$
and $\theta=1$ on $G_1$. Using $\theta^2 w_n\in H^1_c(G_2)$ as a test function we obtain
$$\int_{G_2}w_n^{p+1}\theta^2=
\int_{G_2}\theta^2\nabla w_n\cdot a\cdot\nabla w_n\:dy+
2\int_{G_2}\theta w_n\nabla w_n\cdot a\cdot\nabla\theta\:dy.$$
Thus, by standard computations
\begin{eqnarray*}
\frac{1}{2}\int_{G_2}\theta^2\nabla w_n\cdot a\cdot\nabla w_n\:dx
&\le&2\int_{G_2}w_n^2\nabla \theta\cdot a\cdot\nabla\theta\:dx+
\int_{G_2}\theta^2 w_n^{p+1}\:dx\\
&\le&2c_1\|(\nabla\theta)^2\|_\infty\int_{G_2}u^2\:dx+\int_{G_2}u^{p+1}\:dx.
\end{eqnarray*}
We conclude that $(w_n)$ is bounded in $H^1(G_1)$.
By the construction $v_\psi\le w_n\le u\in H^1(G_1)$ for all $n\in\N$.
Therefore $(w_n)$ has a subsequence, denoted by $(w_{n_1(k)})_{k\in\N}$,
which converges to a function $w^{(1)}\in H^1(G_1)$ weakly in $H^1(G_1)$,
strongly in $L^2(G_1)$ and almost everywhere in $G_1$.
Hence it is clear that $w^{(1)}$ is a solution to (\ref{*}) in $G_1$
and $v_\psi\le w^{(1)}\le u$.

Now we proceed by the standard diagonal argument (see, e.g., \cite[Theorem 2.10]{Ni}).
At the second step, consider a sequence $(w_{n_1(k)})_{k\in\N}$ in $G_2$ (assuming that $n_1(1)>2$).
In the same way as above we obtain a subsequence $(w_{n_2(k)})_{k\in\N}$
that converges to a function $w^{(2)}\in H^1(G_2)$, which is a solution to (\ref{*}) in $G_2$.
Moreover, $v_\psi\le w^{(2)}\le u$ in $G_2$ and $w^{(2)}=w^{(1)}$ in $G_1$.
Continuing this process, for each fixed $m>2$ we construct a subsequence $(w_{n_m(k)})_{k\in\N}$
(with $n_m(1)>m$) that converges weakly to $w^{(m)}\in H^1(G_m)$ which is a solution
to (\ref{*}) in $G_m$ and such that $v_\psi\le w^{(m)}\le u$ in $G_m$, $w^{(m)}=w^{(m-1)}$ in $G_{m-1}$.

By a diagonal process $(w_{n_m(m)})_{m\in\N}$ is a subsequence of $(w_{n_m(k)})_{k\in\N}$
for every $m\in\N$. Thus for each fixed $k\in\N$ the sequence $(w_{n_m(m)})$ converges weakly
to $w^{(k)}$ in $H^1(G_k)$. Let $w_\ast$ be the weak limit of $(w_{n_m(m)})$
in $H^1_{loc}(\C_\Omega^r)$. Then $w_\ast$ is a solution of (\ref{*}) in $\C_\Omega^r$
such that $v_\psi\le w_\ast\le u$ in $\C_\Omega^r$.
\qed

\begin{remark}
The constructed solution $w_\ast$ is actually locally H\"older continuous.
Indeed, since $p\le\frac{N}{N-2}<\frac{N+2}{N-2}$ we conclude by the Brezis--Kato estimate
(see, e.g. \cite[Lemma B.3]{Struwe})
that $w_\ast\in L^s_{loc}(\C_\Omega^r)$ for any $s<\infty$.
Then $-\grad{a}w_\ast=w_\ast^p\in L^s_{loc}(\C_\Omega^r)$ for any $s<\infty$.
Hence the standard elliptic estimates imply that $w_\ast\in C^{0,\gamma}_{loc}(\C_\Omega^r)$.
\end{remark}

\section{Proof of Theorems \ref{Bandle} and \ref{domain}}

In this section we study positive supersolutions at infinity to the model equation
\begin{equation}\label{nonlin}
-\Delta u=u^p\quad\mbox{in }\:\C_\Omega,
\end{equation}
where $p>1$ and $\Omega$ is a subdomain of $S^{N-1}$.
Recall, that $\lambda_1$ denotes the principal eigenvalue of the
Dirichlet Laplace--Beltrami operator $-\Delta_\omega$ in $\Omega$
and $\alpha_-$ stands for the negative root of the equation
$\alpha(\alpha+N-2)=\lambda_1$.

Existence of positive supersoltuions to (\ref{nonlin}) with
$p>p^\ast(id,\C_\Omega)=1-2/\alpha_-$ can be easily verified.
Namely, by direct computations one can find supersolutions of the form
$u=cr^{2/(1-p)}\phi$, where $\phi>0$ is the principal eigenfunction of $-\Delta_\omega$ on $\Omega$
(see also \cite{Bandle-Levine,Bandle} for a direct proof of the existence of positive solutions).
We are going to prove nonexistence of positive supersolutions to (\ref{nonlin}) in $\C_\Omega^\rho$
for $p\in(1,1-2/\alpha_-)$. Notice that if $u>0$ is a solution to (\ref{nonlin}) in $\C_\Omega^\rho$
then, by the scaling properties of the Laplacian,
$\rho^{2}u(x/\rho)$ is a solution to (\ref{nonlin}) in $\C_\Omega^1$.
So in what follows we fix $\rho=1$.

\paragraph{Minimal solution estimate.}

Here we derive the sharp asymptotic at infinity of the minimal solutions to the equation
\begin{equation}\label{min-Delta}
-\Delta u-\frac{V(\omega)}{|x|^2}u=0\quad\mbox{in }\:\C_\Omega
\end{equation}
with $V\in L^\infty(\Omega)$.
Let $-\Delta_\omega$ be the Dirichlet Laplace--Beltrami operator in $L^2(\Omega)$
and $0\le V\in L^\infty(\Omega)$.
Let $(\tilde\lambda_k)_{k\in\N}$
be the sequence of eigenvalues of the operator $-\Delta_\omega-V$,
such that $\tilde\lambda_1<\tilde\lambda_2\le\tilde\lambda_3\le\dots$.
By $(\tilde\phi_k)_{k\in\N}$ we denote the corresponding orthogonal basis
of normalized eigenfunctions in $L^2(\Omega)$, with $\tilde\phi_1>0$.

From now on we assume that $\tilde\lambda_1>-(N-2)^2/4$.
Then the roots of the quadratic equation $\alpha(\alpha+N-2)=\tilde\lambda_k$ are real for each $k\in\N$.
By $\tilde\alpha_k$ we denote the smallest root of the equation, i.e.,
$$\tilde\alpha_k:=-\frac{N-2}{2}-\sqrt{\frac{(N-2)^2}{2}+\tilde\lambda_k}.$$
Notice that since $\tilde\lambda_1>-(N-2)^2/4$
it follows from the Hardy inequality (\ref{Hardy})
that the potential $V(\omega)|x|^{-2}$ is form bounded.
Hence for a given $\psi\in C^\infty_c(\Omega)$ the minimal positive solution $v_\psi$
to (\ref{min-Delta}) in $\C_\Omega^1$ can be constructed as in (\ref{psi-tilde}), (\ref{e-psi}).
We are going to show that $v_\psi$ can be represented as a series of $r^{\tilde\alpha_k}\tilde\phi_k$.

\begin{lemma}
Let $v_\psi$ be the minimal positive solution to (\ref{min-Delta}) in $\C_\Omega^1$. Then
\begin{equation}\label{min-sum}
v_\psi(x)=\sum_{k=1}^\infty\psi_k r^{\tilde\alpha_k}\tilde\phi_k(\omega),\qquad\mbox{where}\quad
\psi_k=\int_\Omega\psi(\omega)\tilde\phi_k(\omega)\:d\omega.
\end{equation}
\end{lemma}

\proof
Set $v_k(x):=r^{\tilde\alpha_k}\tilde\phi_k(\omega)$.
Then a direct computation gives that
$$-\Delta v_k-\frac{V(\omega)}{|x|^2} v_k =0\quad\mbox{in }\:\C_\Omega^1.$$
Recall that $\nabla=\nu\frac{\partial}{\partial r}+\frac{1}{r}\nabla_\omega$
where $\nu=\frac{x}{|x|}\in\R^N$.
Since $\|\nabla_\omega\tilde\phi_k\|_2^2=\tilde\lambda_k$ and $\nu\cdot\nabla_\omega=0$ we obtain
\begin{eqnarray*}
\|\nabla v_k\|_{L^2}^2&=&
\int_1^\infty\int_\Omega\left(|\frac{\partial}{\partial r}r^{\tilde\alpha_k}\tilde\phi_k(\omega)|^2+
\frac{|r^{\tilde\alpha_k}\nabla_\omega\tilde\phi_k(\omega)|^2}{r^2}\right)r^{N-1}\:d\omega dr\\
&=&\int_1^\infty r^{2\tilde\alpha_k+N-3}(\tilde\alpha_k^2+\tilde\lambda_k)\:dr
=\frac{\tilde\alpha_k^2+\tilde\lambda_k}{2-N-2\tilde\alpha_k}=|\tilde\alpha_k|;
\end{eqnarray*}
$$\int_{\C_\Omega^\rho}\nabla v_k\cdot\nabla v_n\:dx=0\quad\mbox{for any }\:k\neq n.$$
Now it is straightforward that $v_k-{\tilde\phi_{k}\theta_1}\in D^1_0(\C_\Omega^1)$,
so $v_k$ solves the problem
$$-\Delta v-\frac{V(\omega)}{|x|^2}v=0,\qquad
v-{\tilde\phi_{k}\theta_1}\in D^1_0(\C_\Omega^1).$$
Therefore
$$\|\nabla v_\psi\|_2^2=\sum_{k=1}^\infty
\psi_k^2\|\nabla v_k\|_2^2\le c\sum_{k=1}^\infty\psi_k^2\sqrt{\tilde\lambda_k}\le c\|\psi\|_2\|\nabla\psi\|_2.$$
Hence $v_\psi-\psi\theta_1\in D^1_0(\C_\Omega^1)$, so $v_\psi$ solves the problem
$$-\Delta v-\frac{V(\omega)}{|x|^2}v=0,\qquad
v-{\psi\theta_1}\in D^1_0(\C_\Omega^1).$$
By the uniqueness we conclude that $v_\psi$ defined by (\ref{min-sum})
coincides with the minimal solution $v_\psi$ as constructed in (\ref{psi-tilde}), (\ref{e-psi}).
\qed

\begin{lemma}\label{alpha-sharp}
Let $v_\psi>0$ be a minimal solution to (\ref{min-Delta}) in $\C_\Omega^1$.
Then for any $\Omega^\prime\Subset\Omega$ and $\rho>1$ there exists $c=c(\Omega^\prime,\rho)>0$ such that
\begin{equation}\label{Delta-sharp}
v_\psi(x)\ge c r^{\tilde\alpha_1}\quad\mbox{in }\:\C_{\Omega^\prime}^{\rho}.
\end{equation}
\end{lemma}

\proof
By (\ref{min-sum}) one can represent $v_\psi$ as
$v_\psi(x)=\psi_1 r^{\tilde\alpha_1}\tilde\phi_1(\omega)+w(x)$, where
$$w(x)=\sum_{k=2}^\infty\tilde\psi_k r^{\tilde\alpha_k}\tilde\phi_k(\omega).$$
Notice that $w(x)$ satisfies
$$-\Delta w-\frac{V(\omega)}{|x|^2}w=0\quad\mbox{in }\:\C_\Omega^1.$$
Thus by the standard elliptic estimate (see, e.g. \cite[Theorem 8.15]{Gilbarg})
for any $\Omega^\prime\Subset\Omega$ and $\rho>1$ one has
$$\sup_{\C_{\Omega^\prime}^{(\rho,2\rho)}}|w|^2\le c\rho^{-N}\int_{\C_{\Omega^\prime}^{(\rho,2\rho)}}|w|^2\;dx,$$
where the constant $c>0$ does not depend on $\rho$. Therefore
\begin{eqnarray*}
\sup_{\C_{\Omega^\prime}^{(\rho,2\rho)}}|w|^2 &\le&
c\rho^{-N}\int_\rho^{2\rho}r^{N-1}\int_{\Omega}|w|^2\:d\omega\,dr
\le c\rho^{-N}\int_\rho^{2\rho}r^{N-1}\sum_{k=2}^\infty \psi_k^2 r^{2\tilde\alpha_k}\:dr\\
&\le&
c\int_\rho^{2\rho}r^{2\tilde\alpha_2-1}\:dr\|\psi-\psi_1\|_2^2
\le c_1 \rho^{2\tilde\alpha_2}.
\end{eqnarray*}
So we conclude that
$$v_\psi(x)\ge\psi_1 r^{\tilde\alpha_1}\phi_1(\omega)-cr^{\tilde\alpha_2}\quad\mbox{in }\:\C_{\Omega^\prime}^\rho,$$
Since $\tilde\alpha_2<\tilde\alpha_1<0$ this implies (\ref{Delta-sharp}).
\qed

\paragraph{Proof of Theorem \ref{Bandle}.}
We distinguish the subcritical and critical cases.

\noindent
{\it Subcritical case $1<p<1-2/\alpha_-$}.
Assume that $u>0$ is a supersolution to (\ref{nonlin}) in $\C_\Omega^r$ for some $r\in(0,1)$.
Then $u>0$ is a supersolution to
\begin{equation}\label{min-sub}
-\Delta u=0\quad\mbox{in }\:\C_\Omega^r.
\end{equation}
By Lemma \ref{minimal} we conclude that $u>cv_\psi$ in $\C_\Omega^1$,
where $v_\psi>0$ is a minimal solution to (\ref{min-sub}).
Then by Lemma \ref{alpha-sharp} for a subdomain $\Omega^{\prime}\Subset\Omega$ one has
\begin{equation}\label{alpha-1}
v_\psi\ge c|x|^{\alpha_-}\quad\mbox{in }\:\C_{\Omega^{\prime}}^1.
\end{equation}
So $u>0$ is a supersolution to
\begin{equation}\label{min-V}
-\Delta u-Wu=0\quad\mbox{in }\:\C_{\Omega}^1,
\end{equation}
where $W(x):=u^{p-1}(x)$ satisfies
$$W(x)\ge c^{p-1}|x|^{\alpha_-(p-1)}\quad\mbox{in }\:\C_{\Omega^\prime}^1,$$
with $\alpha_-(p-1)>-2$. Now Lemma \ref{lin-nonexist} leads to a contradiction.
\medskip

\noindent
{\it Critical case $p=1-2/\alpha_-$.}
Let $u>0$ be a supersolution to (\ref{nonlin}) in $\C_\Omega^r$
with the critical exponent $p^\ast=1-2/\alpha_-$.
Then arguing as in the previous case we conclude that
$u$ is a supersolution to (\ref{min-V}) with
$W(x):=u^{p^\ast-1}(x)$ satisfying
$$W(x)\ge\frac{c^{p^\ast-1}}{|x|^2}\quad\mbox{in }\:\C_{\Omega^{\prime}}^1$$
on a subdomain $\Omega^\prime\Subset\Omega$.
Let $\chi_{\Omega^\prime}(\omega)$ be the characteristic function of $\Omega^\prime$.
Then $u$ is a supersolution in $\C_\Omega^1$ to the equation
\begin{equation}\label{min-chi}
-\Delta v-\frac{\epsilon\chi_{\Omega^\prime}(\omega)}{|x|^2}v=0\quad\mbox{in }\:\C_\Omega,
\end{equation}
for any $\epsilon\in[0,c^{p^\ast-1}]$.
By the variational characterization of the principal Dirichlet eigenvalue
one can fix $\epsilon>0$ small enough in such a way that
$\tilde\lambda_1=\tilde\lambda_1(-\Delta_\omega\!-\!\epsilon\chi_{\Omega^\prime},\Omega)>-(N-2)^2/4$.
Let $w_\psi$ be a minimal positive solution in $\C_\Omega^1$ to (\ref{min-chi})
with such fixed $\epsilon$.
Applying Lemma \ref{alpha-sharp} to (\ref{min-chi})
we conclude that for a subdomain $\Omega^{\prime\prime}\Subset\Omega$ one has
$$w_\psi\ge c|x|^{\tilde\alpha_1}
\quad\mbox{in }\:\C_{\Omega^{\prime\prime}}^1,$$
where $\tilde\alpha_1>\alpha_-$.
So $u>0$ is a supersolution to
\begin{equation}\label{super-V}
-\Delta u-Wu=0\quad\mbox{in }\:\C_{\Omega}^1,
\end{equation}
where $W(x):=u^{p^\ast-1}(x)$ satisfies
$$W(x)\ge c^{p^\ast-1}|x|^{\tilde\alpha_1(p^\ast-1)}\quad\mbox{in }\:\C_{\Omega^{\prime\prime}}^1$$
with $\tilde\alpha_1(p^\ast-1)>-2$. This contradicts to Lemma \ref{lin-nonexist}.
\qed

\begin{remark}
Strictly speaking, in the above proof the subcritical case $1<p<1-2/\alpha_-$
is redundant, due to Proposition \ref{properties} (ii).
\end{remark}

Let $\Omega\subset S^{N-1}$ be a domain
such that $\lambda_1=\lambda_1(\Omega)>0$.
Define the operator $L_d$ by
\begin{equation}\label{Ld}
L_d=-\frac{\partial^2}{\partial r^2}
-\frac{N-1}{r}\frac{\partial}{\partial r}-\frac{d(r)}{\lambda_1}\frac{1}{r^2}\Delta_\omega,
\end{equation}
where $d(r)$ is measurable and squeezed between two positive constants.
Then $L_d$ is a divergence type uniformly elliptic operator $-\grad{a_d}$ (see, e.g., \cite{Vogt}).

\paragraph{Proof of Theorem \ref{domain}.}
Consider the operator $L_d$ where $d(r)\equiv \alpha(\alpha+N-2)$ with $\alpha<2-N$.
Following the lines of the proof of Theorem \ref{Bandle}
we conclude that $p^\ast(a_d,\C_\Omega)=1-2/\alpha$.
Clearly for any given $p\in(1,\frac{N}{N-2})$, one can choose $\alpha$ such
that $p^\ast(a_d,\C_\Omega)=p$.
\qed

\begin{remark}\label{R-log}
In the above proof equation (\ref{*}) has no positive
supersolutions at infinity in $\C_\Omega$ in the critical case
$p=p^\ast(a_d,\C_\Omega)$. Next we give an example of equation
(\ref{*}) with a positive supersolution at infinity in the
critical case.

Let $\Omega\subset S^{N-1}$ be smooth and $L_{\tilde d}$ be as in (\ref{Ld}) with
$$\tilde d(r)=\alpha(\alpha+N-2)+\frac{2-N-2\alpha}{\log(r)}+\frac{2}{\log^2(r)},$$
where $\alpha<2-N$.
For large enough $R\gg 1$ the operator $L_{\tilde d}=-\grad{a_{\tilde d}}$ is uniformly elliptic on $\C_\Omega^R$.
Let $\phi_1>0$ be the principal Dirichlet eigenfunction of $-\Delta_\omega$,
corresponding to $\lambda_1$.
Direct computation shows that the function
$$v_{\phi_1}:=\frac{r^\alpha}{\log(r)}\phi_1$$
is a solution to the equation
\begin{equation}\label{Ld0}
L_{\tilde d} v=0\quad\mbox{in }\:\C_\Omega^R.
\end{equation}
Since $\Omega$ is smooth, the Hopf lemma implies that $v_{\phi_1}$
is a minimal solution to (\ref{Ld0}), in the sense of Lemma \ref{minimal},
that is, for any supersolution $u>0$ to (\ref{Ld0}) in $\C_\Omega^R$ and for any $\rho>R$,
there exists $c>0$ such that $u\ge cv_{\phi_1}$ in $\C_\Omega^\rho$.
Following the lines of the proof of Theorem \ref{Bandle}, subcritical case,
we conclude that $p^\ast(a_d,\C_\Omega)=1-2/\alpha$.
On the other hand, one can readily verify that $u=r^\alpha\phi_1$
is a positive supersolution to (\ref{*}) in the critical case $p=1-2/\alpha$.

Note that the value of the critical exponent for $L_{\tilde d}$ is the same as $L_d$
due to the fact that $\lim\limits_{r\to\infty}({\tilde d}(r)-d(r))=0$.
However the rate of convergence is not sufficient
to guarantee the equivalence of the corresponding minimal solutions
(see, e.g. \cite{Ancona,Pinchover-Green} for the related estimates of Green's functions).
This explains the nature of the different behaviour of the nonlinear equations (\ref{*})
at the critical value of $p$.
\end{remark}

\section{Proof of Theorem \ref{matrix}}

First we show that for any domain $\Omega\subseteq S^{N-1}$ one has $p^\ast(a,\C_\Omega)>1$.
Then we prove the second part of theorem \ref{matrix},
saying that if the complement of $\Omega$ has nonempty interior then $p^\ast(a,\C_\Omega)<\frac{N}{N-2}$.
We start with establishing a lower bound on positive solutions of the equation
\begin{equation}\label{min-a}
-\grad{a}v=0\quad\mbox{in }\:\C_\Omega.
\end{equation}

\begin{lemma}\label{alpha}
Let $\Omega\subseteq S^{N-1}$ be a domain and $\Omega^\prime\Subset\Omega$.
Then there exists $\alpha=\alpha(\Omega^\prime)\le 2-N$ such that for any $\rho>0$
any positive solution $v$ to equation (\ref{min-a}) in $\C_\Omega^\rho$
has a polynomial lower bound
\begin{equation}\label{lower}
v\ge c|x|^\alpha\quad\mbox{in }\:\C_{\Omega^\prime}^{2\rho}.
\end{equation}
\end{lemma}

\proof
Set $a=3/4$, $b=7/4$. Let $r\ge 2\rho$ and $m_r=\inf_{\C_{\Omega^\prime}^{(ra,rb)}}v$.
By the strong Harnack inequality $v$ satisfies
$$\inf_{\C_{\Omega^\prime}^{(ra,rb)}} v\ge C_S
\sup_{\C_{\Omega^\prime}^{(ra,rb)}} v,$$
with the constant $C_S\in(0,1)$ dependent on $\Omega^\prime$ and not on $r$,
as a simple scaling argument shows.
Then
\begin{equation}\label{chain}
m_r\le\sup_{\C_{\Omega^\prime}^{(ra,rb)}}v
\le\sup_{\C_{\Omega^\prime}^{(ra,2rb)}}v\le
C_S^{-1}\inf_{\C_{\Omega^\prime}^{(ra,2rb)}}v\le
C_S^{-1}\inf_{\C_{\Omega^\prime}^{(2ra,2rb)}}v=C_S^{-1}m_{2r}.
\end{equation}
Let $r_n=2^n\rho$ and $n\in\N$.
Iterating (\ref{chain}) we obtain $m_{r_n}\ge C_S^{n-1}m_{2\rho}$.
Choosing $n$ such that $ar_n\le|x|<2ar_n$ one can see that
$$v\ge c|x|^\alpha\quad\mbox{in }\:\C_{\Omega^\prime}^{2a\rho},$$
where $\alpha=\log_2 C_S$ and $c=c(\rho)=(a\rho)^{-\alpha}C_S^{-1}m_{2\rho}$.
Taking into account (\ref{upper-green}) we conclude that $\alpha\le 2-N$.
\qed

\begin{remark}
A similar argument was used before by Pinchover \cite[Lemma 6.5]{Pinchover}.
Observe that in the same way one can get a rough polynomial upper bound
on positive solutions of (\ref{min-a}).
\end{remark}

The lower bound bound (\ref{lower}) allows us
to prove nonexistence of positive solutions to (\ref{*})
exactly by the same argument as was used in the proof of  Theorem \ref{Bandle}
in the subcritical case.

\begin{proposition}\label{nonexist}
Let $\Omega\subseteq S^{N-1}$ be a domain.
Then $p^\ast(a,\C_\Omega)\ge 1-2/\alpha$ where $\alpha\le 2-N$ is from the lower bound (\ref{lower}).
\end{proposition}

\proof
Assume that $u\ge 0$ is a supersolution to (\ref{*}) in $\C_\Omega^\rho$ with exponent $p<1-1/\alpha$.
By Lemma \ref{comparison} and (\ref{lower}) we conclude that
for any subdomain $\Omega^\prime\subset\Omega$ there exists $c=c(\Omega)>0$ such that
$$u\ge c|x|^\alpha\quad\mbox{in }\:\C_{\Omega^\prime}^{2\rho+2}.$$
Therefore $u$ is a supersolution to
$$-\grad{a}u=Vu\quad\mbox{in }\:\C_{\Omega^\prime}^{2\rho+2},$$
where $V(x):=u^{p-1}(x)$ satisfies the inequality
$$V(x)\ge c^\prime|x|^{\alpha(p-1)}\quad\mbox{in }\:\C_{\Omega^{\prime}}^{2\rho+2}$$
with $\alpha(p-1)>-2$.
Then Lemma \ref{lin-nonexist} implies that $u\equiv 0$ in $\C_\Omega^\rho$.
Since $\alpha>0$ does not depend on $\rho$, we conclude that $p^\ast(a,\C_\Omega)\ge 1-1/\alpha$.
\qed

Our next step is to obtain a polynomial upper bound
on the minimal positive solutions to the equation
$$-\grad{a}v-Vv=0\quad\mbox{in }\:\C_\Omega,$$
with a special potential $V$ which will be specified later.
In order to do this we need the notion of a Green bounded potential.
Let $\Gamma_a(x,y)$ be the fundamental solution to
$$-\grad{a}v=0\quad\mbox{in }\:\R^N.$$
We say that a potential $0\le V\in L^1_{loc}(\R^N)$ is Green bounded and write $V\in GB$ if
$$\|V\|_{GB,a}:=\sup_{x\in\R^N}\int_{\R^N}\Gamma_a(x,y)V(y)dy<\infty,$$
which is equivalent up to a constant factor to the condition
$\sup_{x\in\R^N}\int_{\R^N}|x-y|^{2-N}|V(y)|dy<\infty$,
but we will use below the numerical value of $\|V\|_{GB,a}$.
One can see, e.g. by the Stein interpolation theorem, that
if $V\in GB$ then $V$ is form bounded in the sense of (\ref{form-bound}).
We will use the following important property of Green bounded potentials,
which was proved in \cite{Hansen}, see also \cite{Pinchover-Green}.
\begin{lemma}\label{Green-bounded}
Let $V\in GB$ and $\|V\|_{GB,a}<1$. Then there exists a solution $w>0$ to the equation
\begin{equation}\label{quasi-const}
-\grad{a}w-Vw=0\quad\mbox{in }\:\R^N,
\end{equation}
such that $0<c<w<c^{-1}$ in $\R^N$.
\end{lemma}

Using this result we first prove the required upper bound in the case of the
"half--space" cone $\C_+=\{x_N>0\}$ with the cross--section $S_+=\{|x|=1,x_N>0\}$.
For a given uniformly elliptic matrix $a$ and a potential $V$ defined on $\C_+$
we denote by $\bar{a}$ and $\bar{V}$ the extensions of $a$ and $V$ from $\C_+$ to $\R^N$ by reflection,
so that $\bar{a}(\cdot,-x_N)=a(\cdot,x_N)$ and $\bar{V}(\cdot,-x_N)=\bar{V}(\cdot,x_N)$.
Thus the matrix $\bar{a}$ is uniformly elliptic on $\R^N$ with the same ellipticity constant as $a$.

\begin{lemma}\label{lemma-half-space}
Let $0\le V\in L^1_{loc}(\C_+)$ be a potential such that $\|\bar{V}\|_{GB,\bar{a}}<1$.
Let $v_\psi>0$ be a minimal positive solution in $\C_+^1$ to the equation
\begin{equation*}
-\grad{a}v-Vv=0\quad\mbox{in }\:\C_+.
\end{equation*}
Then there exists $\gamma\in(0,1)$ such that
\begin{equation}\label{half-space}
0<v_\psi\le c|x|^{2-N-\gamma}\quad\mbox{in }\:\C_+^1.
\end{equation}
\end{lemma}

\proof
Let $\bar{v}$ denote the extension of $v_\psi$ from $\C_+^1$ to $B_1^c$ by reflection,
that is $\bar{v}(\cdot,x_N)=-v_\psi(\cdot,-x_N)$.
Thus $\bar{v}(x)$ is a solution to the equation
\begin{equation*}
-\grad{\bar a}\bar{v}-\bar{V}\bar{v}=0\quad\mbox{in }\:B_1^c.
\end{equation*}
Let $w$ be a solution to (\ref{quasi-const}) given by Lemma \ref{Green-bounded}.
One can check by direct computation (see \cite[Lemma 3.4]{KLS}),
that $v_1:=\bar{v}/w$ is a solution to the equation
\begin{equation}\label{bar}
-\grad{(w^2\bar{a})}v_1=0\quad\mbox{in }\:B_1^c,
\end{equation}
where the matrix $w^2\bar{a}$ is clearly uniformly elliptic.
Let $\Gamma(x):=\Gamma_{w^2\bar{a}}(x,0)$ be the fundamental solution
to the equation $-\grad{(w^2\bar{a})}u=0$ in $\R^N$.
By the classical estimate \cite{LSW} one has
\begin{equation}\label{w-1}
c_1|x|^{2-N}\le \Gamma(x)\le c_2|x|^{2-N}\quad\mbox{in }\:B_1^c.
\end{equation}
Applying Lemma \ref{comparison} to $v_1$ and $\Gamma$ on $C_+^1$
and by the construction of $v_1$ we conclude that
\begin{equation}\label{v-1}
|v_1(x)|\le c_3\Gamma(x) \quad\mbox{on }\:B_1^c.
\end{equation}
Applying the Kelvin transformation $y=y(x)=x/|x|^2$ and $x=x(y)=y/|y|^2$
to (\ref{bar}) we see that the function $\tilde{v}_1(y)=v_1(x(y))/\Gamma(x(y)), \, \tilde{v}\in L^\infty(B_1)$,
solves the equation
\begin{equation*}
-\grad{\tilde a}\tilde{v}_1=0\quad\mbox{in }\:B_1,
\end{equation*}
where the matrix $\tilde{a}(y)$ is uniformly elliptic on $B_1$.
It follows that $\tilde{v}_1\in H^1_{loc}(B_1)$ (see, e.g., \cite{Serrin}).
Then, by the De Giorgi -- Nash regularity result \cite{Gilbarg},
$\tilde{v}_1\in C^{0,\gamma}(B_1)$ for some $\gamma\in(0,1)$.
Notice that
\begin{equation*}\label{boundary}
\tilde{v}_1(y)=0\quad\mbox{in}\quad\{y\in B_1,\:y_N=0\}
\end{equation*}
by the construction. Therefore
$\tilde{v}_1(0)=0$, hence
$$|\tilde{v}_1(y)|\le c|y|^\gamma\quad\mbox{in }\:B_1.$$
We conclude that
$$|\bar{v}|\le c_3|\bar{v}_1(x)|\le c_4|x|^{2-N-\gamma}
\quad\mbox{in }\:B_1^c,$$
as required.
\qed

\begin{lemma}\label{beta}
Let $\Omega\subset S^{N-1}$ be a domain such that $S^{N-1}\setminus\Omega$ has nonempty interior.
Let
$$W_\epsilon(x):=\frac{\epsilon}{|x|^2\log^2(|x|+2)}.$$
Then there exists $\epsilon>0$ and $\beta=\beta(\epsilon)<2-N$
such that any minimal positive solution $v_\psi$ in $\C_\Omega^1$ to the equation
\begin{equation*}
-\grad{a}v-W_\epsilon v=0\quad\mbox{in }\:\C_\Omega
\end{equation*}
has the polynomial upper bound
\begin{equation}\label{upper}
v_\psi\le c|x|^\beta\quad\mbox{in }\:\C_\Omega^1.
\end{equation}
\end{lemma}

\proof
If $\C_\Omega\subseteq\C_+$ then (\ref{upper})
follows from (\ref{half-space}) by Lemma \ref{comparison}.
We shall consider the case $\C_\Omega\not\subseteq\C_+$.

Without loss of generality we can assume that $(0,\dots,0,-1)\not\in\Omega$.
Set $\hat x=(x_1,\dots,x_{N-1})$ and $\sigma=\inf\{|\hat x|:x\in\Omega,x_N<0\}$.
Let
$D_\sigma=\{x\in S^{N-1}:|\hat x|\le\sigma,x_N<0\}$ and $\hat D_\sigma:=S^{N-1}\setminus D_\sigma$.
Then $\C_\Omega\subseteq\C_{\hat D_\sigma}$.
Extend the matrix $a$ by $id$ from $\C_\Omega$ to $\C_{\hat D_\sigma}$.
Let $w_\psi$ be a minimal positive solution in $\C_{\hat D_\sigma}^1$ to the equation
\begin{equation*}
-\grad{a}w-W_\epsilon w=0\quad\mbox{in }\:\C_{\hat D_\sigma}
\end{equation*}
To complete the proof we need only to show that $w_\psi$ satisfies (\ref{upper}) in $\C_{\hat D_\sigma}^1$.
Then the same bound on minimal solutions in $\C_\Omega^1$ follows from Lemma \ref{comparison}.

Consider the transformation
$$y=y(x)=\left(x_1,\dots,x_{N-1},x_N+k|\hat x|\right),$$
where $k=\sqrt{\sigma^{-2}-1}$.
Then $y:\C_{\hat D_\sigma}\to\C_+$ is a bijection,
the Jacobian of $y(x)$ is nondegenerate and has the determinant equal to $1$ everywhere.
Moreover, $|x|\le|y(x)|\le \kappa|x|$ for all $x\in\C_{\hat D_\sigma}$, where $\kappa=\sqrt{2+k^2}$.
Therefore $\hat{w}(y):=w_\psi(x(y))$ solves the equation
$$-\grad{\hat a}\hat{w}-\hat W_\epsilon \hat{w}=0\quad
\mbox{in $\C_+^\kappa$,}$$
with the uniformly elliptic matrix $\hat{a}(y):=a(x(y))$
and $\hat W_\epsilon(y):=W_\epsilon(x(y))$.
One can easily check by direct computation that
$\bar{\hat W}_\epsilon\in GB$.
Fix $\epsilon>0$ such that $\|\bar{\hat W}_\epsilon\|_{GB,\bar{\hat{a}}}<1$.
Then by Lemma \ref{lemma-half-space} we conclude that $\hat{w}(y)$ satisfies (\ref{half-space}).
Therefore $w_\psi(x)$ obeys (\ref{upper}) with $\beta:=2-N-\gamma$ as required.
\qed

\begin{proposition}\label{exist}
Let $\Omega\subset S^{N-1}$ be a domain such that $S^{N-1}\setminus\Omega$ has nonempty interior.
Then $p^\ast(a,\C_\Omega)\le 1-2/\beta$, where $\beta<2-N$ is from the upper bound (\ref{upper}).
\end{proposition}

\proof
Fix $p>p_0=1-2/\beta$ and set $\delta=p-p_0$.
Let $w_\psi>0$ be a minimal positive solution in $\C_\Omega^1$ to
$$-\grad{a}w-W_\epsilon w=0\quad\mbox{in }\:\C_\Omega$$
where $\epsilon>0$ is from Lemma \ref{beta}.
Then by (\ref{upper}) for some $\bar\tau=\bar\tau(\delta)>0$ small enough the function
$\bar\tau w_\psi$ satisfies
$$(\bar\tau w_\psi)^{p-1}\le \bar\tau^{p-1}(c|x|^\beta)^{p-1}\le
\frac{\bar\tau^{p-1}c_1}{|x|^{2+\delta|\beta|}}\le
\frac{\epsilon}{|x|^2\log^2(|x|+2)}=
W_\epsilon(x)\quad\mbox{in }\C_\Omega^1.$$
Therefore
$$-\grad{a}(\bar\tau w_\psi)=W_\epsilon(\bar\tau w_\psi)\ge
(\bar\tau w_\psi)^{p-1}(\bar\tau w_\psi)=(\bar\tau w_\psi)^p
\quad\mbox{in }\C_\Omega^1,$$
that is $\bar\tau w_\psi>0$ is a supersolution to (\ref{*}) in $\C_\Omega^1$.
\qed

\paragraph{Concluding remarks.}
The proofs of Propositions \ref{nonexist} and \ref{exist}
rely only on the polynomial lower and upper bounds (\ref{lower}) and (\ref{upper}).
Namely, given $\alpha\le\beta<2-N$ in (\ref{lower}) and (\ref{upper})
we conclude that
$$1-\frac{2}{\alpha}\le p^\ast(a,\C_\Omega)\le 1-\frac{2}{\beta}.$$
By the next example we show
that the (optimal) constants $\alpha$ and $\beta$
might be actually different.

Let $\Omega\subset S^{N-1}$ be smooth and $L_d$ be as in (\ref{Ld}) with
$$d(r)=A(r)(A(r)+N-2)+R(r),$$
where
$$A(r)=\gamma+\delta[\sin(k\log\log(r))+k\cos(k\log\log(r))],$$
$$R(r)=k\delta[\cos(k\log\log(r))-k\sin(k\log\log(r))]\log^{-1}(r),$$
$\gamma<2-N$, $\delta>0$ and $k>0$ such that $\gamma+\delta\sqrt{k^2+1}<2-N$.
Thus for large enough $R\gg 1$ the operator $L_d=-\grad{a_d}$ is uniformly elliptic on $\C_\Omega^R$.
Let $\phi_1>0$ be the principal Dirichlet eigenfunction of $-\Delta_\omega$,
corresponding to $\lambda_1$.
Direct computation and the Hopf Lemma show that the function
$$v_{\phi_1}:=r^{\gamma+\delta\sin(k\log\log(r))}\phi_1$$
is a minimal solution to the equation
$L_d v=0$ in $\C_\Omega^R$, in the sense of Lemma \ref{minimal}.
Clearly any $\alpha$ and $\beta$ ($\alpha<\beta<2-N$) could be represented as $\alpha=\gamma-\delta$
and $\beta=\gamma+\delta$ for an appropriate choice of parameters $\gamma$, $\delta$ and $k$.
Therefore one cannot expect a sharp polynomial asymptotic of minimal solutions
to the equation $-\grad{a} v=0$ in cone--like domains without additional restrictions on the matrix $a(x)$.

It is an interesting open problem to determine the value of the critical exponent $p^\ast(a,\C_\Omega)$
in the case of minimal solutions oscillating at infinity between two different polynomials.

\begin{small}
\section*{Acknowledgements}
The research of the first named author was supported by the Institute of Advanced Studies
of the University of Bristol via Benjamin Meaker Fellowship.
It is a pleasure to thank the university for support and hospitality.
The authors are grateful to Zeev Sobol for useful comments
and to Yehuda Pinchover for interesting discussions.
\end{small}

\begin{small}

\end{small}

\end{document}